\documentclass[12pt]{article}
\usepackage{graphicx}
\usepackage{hyperref}
\usepackage{amsmath}
\usepackage{amscd,amssymb}

\newcommand{\cF}{\mathcal{F}}
\newcommand{\cH}{\mathcal{H}}

\newcommand{\nP}{\textsf{P}}

\def\cR{{\mathcal R}}
\def\cS{{\mathcal S}}

\def\cD{{\mathcal D}}

\def\cH{{\mathcal H}}

\def\cL{{\mathcal L}}

\def\cF{{\mathcal F}}

\def\cO{{\mathcal O}}

\def\mD{{\mathfrak D}}

\def\mD{{\mathfrak D}}

\def\mS{{\mathfrak S}}

\newtheorem{theorem}{Theorem}[section]

\newtheorem{remark}{Remark}[section]

\newcommand{\beq}{\begin{eqnarray}}
\newcommand{\eeq}{\end{eqnarray}}
\numberwithin{equation}{section}

\begin{document}


\begin{center}
{\large\bf The Duistermaat-Heckman Formula with Application to Circle Actions and
Poincar\'{e} $q$-Polynomials in Twisted Equivariant K-Theory}
\end{center}

\vspace{0.1in}

\begin{center}
{\large
A. A. Bytsenko $^{(a)}$
\footnote{E-mail: aabyts@gmail.com},
M. Chaichian $^{(b)}$
\footnote{E-mail: masud.chaichian@helsinki.fi},
and A. E. Gon\c{c}alves $^{(a)}$}
\footnote{E-mail: aedsongoncalves@gmail.com}
\end{center}

\begin{center}
\vspace{5mm}

$^{(a)}$
{\it
Departamento de F\'{\i}sica, Universidade Estadual de
Londrina\\ Caixa Postal 6001,
Londrina-Paran\'a, Brazil}

\vspace{0.2cm}
$^{(b)}$
{\it
Department of Physics, University of Helsinki\\
P.O. Box 64, 00014 Helsinki, Finland}

\end{center}

\vspace{0.1in}

\begin{abstract}
In this paper we deduce the sketch of proof of the Duistermaat-Heckman formula and investigate
how the known Duistermaat-Heckman result could be specialized to the symplectic structure on the orbit
space. The theorems of localization in equivariant cohomology not only provide us with beautiful mathematical
formulas and stimulate achievements in algorithmic computations, but also promote progress in
theoretical and mathematical physics. We present the elliptic genera and the characteristic
$q$-series for the circle actions and twisted equivariant K-theory, with the case of the
symmetric group of $n$ symbols separately analyzed. We show that the Poincar\'{e} $q$-polynomials
admit presentation in terms of the Patterson-Selberg (or the Ruelle-type) spectral functions.
\end{abstract}

\vspace{0.1in}

\begin{flushleft}
PACS \, 03.70.+k, 11.10.-z (Quantum field theory);  \\
MSC \, 05A30 (q-Calculus and related topics); 19L47 (Equivariant K-theory); \\
19L50 (Twisted K-theory)

\vspace{0.3in}
\end{flushleft}

\newpage

\tableofcontents


\section{Introduction}
\label{Introduction}

One of the aims  of the present  paper is to explain the sketch of the proof of the Duistermaat-Heckman theory
\cite{DH}. Berline and Vergne showed how to derive the Duistermaat-Heckman formula from the
localization theorem for equivariant cohomology (see for example \cite{BGV}, Theorem 7.19).
The theorems of localization in equivariant cohomology not only provide beautiful mathematical formulas
and stimulate achievements in algorithmic computations, but also promote progress in theoretical physics.

One of the remarkable insights of the Duistermaat-Heckman formula is that we can deduce
from it in a consistent way the Itzykson-Zuber formula for arbitrary compact connected Lie group action.
The Itzykson-Zuber integration formula \cite{IZ}, for example, occurs essentially in matrix models
(the Ising model on a random surface) where one considers the coupling of conformal matter to
two-dimensional quantum gravity \cite{Witten92}. This formula also appears in works  on higher-dimensional
lattice gauge theories.

The  specialization of the Duistermaat-Heckman formula leads to the Kirillov integral formula
for irreducible representations of a compact connected Lie group \cite{Semenov}, which has relevance for
the geometric quantization: integration over the matrix groups $G$, which amounts to integration over $G/T$ for
$T$-invariant functions, has well-known importance in  diverse areas  such as integrable models, low-dimensional
gauge theories, quantum gravity and quantum chromodynamics \cite{Blau,Szabo}.

Extension of the Duistermaat-Heckman formula to an infinite-dimensional manifold, the loop space
of smooth maps from the circle $S^1$ to a compact manifold $X$, has been proposed for the first time
by E. Witten. A formal application of the Duistermaat-Heckman theory to the partition function
of $N=1/2$ supersymmetric quantum mechanics reduces to a correct formula for the Dirac index \cite{Atiyah}.
Subsequent arguments supporting  these ideas (with strong mathematical base) were presented in \cite{Bismut,Bismut2}.

Equivariant methods also can be used to define and analyze (twisted) K-theory \cite{Adem,Wang,Freed}.
Applying techniques from the equivariant K-theory, one can construct a Chern character and K-theory Euler class.
This manifests in the fact that the $q$-series elliptic genera can be expressed in terms of $q$-analogs
of the classical special functions, in particular the spectral Patterson-Selberg and Ruelle
functions \cite{Bytsenko11}. In the present work  we would like to call the reader's attention to this connection
by discussing some mathematical aspects associated with it.

{\bf The plan of this paper is as follows:}
In Section \ref{D-H} we begin with discussion of the Pfaffians.
Note that particular attention should be paid  to the choice of sign in the square root
${\rm det}^{1/2}(\xi)$, which depends on the choice of the orientation of a vector field $\xi$.
\footnote{
We will frequently use the notation ${\rm det}^{1/2}(\xi)= {\rm Pf}(\xi)$.
}
We shall investigate the Hamiltonian $G$-spase with  attention to the symplectic structure on
the orbit space, which is a necessary subject  for quantization of the symplectic manifolds. Then we turn
to the Duistermaat-Heckman formula.

In Section \ref{D} we use the language of $\mD$-modules -- sheaves of modules over the sheaf of linear
differential operators. That  allows us to consider some applications of the localization formula and
investigate how the known Riemann-Roch-Hirzebruch integral formula could be generalized to $\mD$-modules.

Section \ref{Circle} is  devoted  to the  application to  the one-dimensional setting of circle actions
on a symplectic manifold.
In Section \ref{N} the developments of the direct analog of K-theory Euler class
make it possible to give a rigorous construction of the Chern character
and the higher elliptic genera of level $N$ in terms of q-series. We show that the $q$-series admit a presentation
in the form of Ruelle-type spectral functions $\cR(s)$. The symmetry properties
(or more general the modular properties) of  appropriate quantities can be analyzed by means of the
functional equations for the spectral functions $\cR(s)$. Results of this section are one of the
novelties of the present work.

Finally, in Section \ref{Twist} we discuss the twisted equivariant K-theory. The Poincar\'{e} polynomials
are  analyzed in Section \ref{Sym}. Using the decomposition formula for the twisted K-theory,
we discover the $q$-series for twisted symmetric products in terms of $\cR(s)$ functions. Further on we turn to
the orbifold symmetric product in Section \ref{Orbifold}. The orbifold Poincar\'{e} $q$-polynomials admit
a presentation in terms of the Ruelle-type spectral functions.

\section{Symplectic structure and the Duistermaat-Heckman formula}
\label{D-H}

{\bf The Pfaffian.}
Suppose that manifold $X$ is oriented, ${\rm dim} X =2n$.
Let $L_p(\xi)$ be a non-singular linear operator on $T_p X: {\rm det} L_p(\xi)\neq 0$, where $\xi$
is a smooth vector field on $X$ and  $T_pX$ is the tangent space of $X$ at the point $p$.
Find an ordered orthonormal basis
${\bf e} = {\bf e}^{(p)}=\{e_j=e_j^{(p)}\}_{j=1}^{2n}$ of
$T_p X $ such that
\begin{equation}
L_p(\xi)e_{2j-1} =  \lambda_je_{2j}\,,\,\,\,\,\,\,\,
L_p(\xi)e_{2j}  =  -\lambda_je_{2j-1},\,\,\,\,\,{\rm for}
\,\,\,\,\, 1\leq j\leq n\,,
\label{L3}
\end{equation}
where each $\lambda_j\in {\mathbb R}/\{0\}$.
Relative to ${\bf e}$, $L_p(\xi)$ is the skew-symmetric matrix.
The Pfaffian ${\rm Pf}_{\bf e}(L_p(\xi))$ of $L_p(\xi)$ associated to $\bf e$ has the
form \cite{BGV}:
$
{\rm Pf}_{\bf e}(L_p(\xi))^2 = {\rm det} L_p(\xi).
$

Assume that $G$ is a compact Lie group which acts
smoothly on $X$ (say on the left). Suppose that the metric $\langle\, ,\, \rangle$ is
$G$-invariant and ${\mathfrak g}$ is the Lie algebra of $G$. For
${\xi}\in {\mathfrak g}$ there is an induced vector field ${\xi}^{*}\in VX$ (= the space of smooth
vector fields on $X$).
A vector field ${\xi}^{*}$ is said to be
non-degenerate if for every zero $p\in X$ of ${\xi}^{*}$ the induced linear map
$L_p({\xi}^{*}): T_p X\rightarrow T_p X$ is non-singular.
${\xi}^{*}$ is a Killing vector field and therefore $L_p({\xi}^{*})$
is skew-symmetric with respect to the inner product structure $\langle \,,\, \rangle_p$
on $T_p(X)$ and the non-singularity of $L_p({\xi}^{*})$ means that one can construct the
square-root
\begin{equation}
\left({\rm det}L_p({\xi}^{*})\right)^{1/2}=(-1)^n{\rm Pf}_{\bf e}(L_p({\xi}^{*}))= \prod_j^n\lambda_j\,.
\label{det3}
\end{equation}
\begin{remark}
The Euler number of an even-dimensional oriented manifold $X$ with Riemannian curvature $R$
is given by the Gauss-Bonnet-Chern formula
$
{\rm Eul}(X) = (2\pi)^{-n/2}\int_X {\rm Pf}(-R)\,.
$
The reader can find the the proof of this basic theorem of differetial geometry,
for example, in the excellent book \cite{BGV}, Section 1.6.
\end{remark}

{\bf symplectic structure on orbit space.}
Suppose $X$ has a symplectic structure $\sigma:\sigma \in \Lambda^2X$ which is a closed two-form
, i.e. $d\sigma=0$, such that for every $p\in X$ the corresponding skew-symmetric form
$
\sigma_p: T_p(X)\oplus T_p(X)\rightarrow {\mathbb R}
$
is non-degenerate. For an integer $j\geq 0$, $\Lambda^jX$ denote the space of smooth complex
differential forms of degree $j$ on smooth manifold $X$. $d\xi: \Lambda^jX\rightarrow \Lambda^{j+1}X$
denotes the exterior differentiation by $\xi$, while $\imath\xi: \Lambda^jX\rightarrow \Lambda^{j-1}X$
denotes the  interior differentiation by $\xi$.
Suppose $X$ is oriented by the Liouville form
$
\omega_{\sigma}=
\frac{1}{n!}\sigma \wedge\cdot\cdot\cdot \wedge\sigma \in
\Lambda^{2n}X/\{0\}\,;
$
suppose also that there is a map
${\mathcal H} : {\mathfrak g}\rightarrow C^{\infty}(X)$
which satisfies
\begin{equation}
\iota ({\xi})\sigma + d{\mathcal H}({\xi})=0,\,\,\,\,\,\,\,\,\,
\forall {\xi}\in {\mathfrak g}\,.
\label{Ham}
\end{equation}
The existence of such a map ${\cH}$ is compatible with the assumption that the action of $G$ on $X$
is Hamiltonian. Let ${\cH} VX$ denote the space of Hamiltonian vector fields on $X$
(actually ${\cH} VX$ is a Lie algebra).

The (left) action of $G$ on $X$ is called symplectic if
${\xi}\in {\cH} VX,\,\,\, \forall {\xi}\in {\mathfrak g}$.
The (left) action of $G$ on $X$ is called
Hamiltonian if it is symplectic and if the Lie algebra homomorphism
$\eta: \,{\mathfrak g}\rightarrow {\mathcal H} VX$ has a lift to
$C^{\infty}(X)$. We note that such a ${\mathcal H}$ will indeed satisfy
condition (\ref{Ham}). The triple $(X,\sigma, {\mathcal H})$ is called a Hamiltonian
$G$-space \cite{Kostant,Woodhouse}.
\begin{theorem} \label{TH1}
(see \cite{BLW}, Theorem 1.)
Assume that $X$ and $G$ are compact
and the Riemannian metric $\langle\, ,\, \rangle$ on $X$ is $G$-invariant; i.e.
each $a\in G$ acts as an isometry of $X$.  Assume now that the induced vector field ${\xi}^{*}$ on $X$ is non-degenerate;
thus the square-root in {\rm (\ref{det3})} is well-defined (and is non-zero)
for $p\in X$ a zero of ${\xi}^{*}$ (i.e. ${\xi}_p^{*}=0$).
Then for any cohomology class $[\tau]\in \Sigma(X,{\xi}^{*},s)$, $s^{-1}=-2\pi i$, one has
\begin{equation}
\int_X[\tau]=(-1)^{n/2}
\!\!\!\!\!
\sum_{\scriptstyle p\in X,
\atop\scriptstyle p=\,{\rm a\,\, zero\,\, of}\,\,{\xi}^{*}}
\frac{p^{*}[\tau]}
{[{\rm det}L_p({\xi}^{*})]^{1/2}}\,.
\label{For1}
\end{equation}
\end{theorem}

{\bf The Duistermaat-Heckman formula.}
We are now in the position to present the Duistermaat-Heckman formula in a form
directly derivable from the Theorem \ref{TH1}.
\begin{theorem}
\label{DHformula}
Suppose $(X,\sigma, {\mathcal H})$ is a Hamiltonian $G$-space, where $G$ and $X$
are compact, and $X$ is oriented by the Liouville form
$
\omega_{\sigma}=
\frac{1}{n!}\sigma \wedge\cdot\cdot\cdot \wedge\sigma \in
\Lambda^{2n}X/\{0\}\,.
$
For $c\in {\mathbb C}$ and for ${\xi}\in {\mathfrak g}$ with ${\xi}^*$ non-degenerate, we have
\begin{equation}
\int_Xe^{c {\mathcal H} ({\xi})}\omega_\sigma = \,\,\, \left(\frac{2\pi}{c}\right)^n
\!\!\sum_{\scriptstyle p\in X
\atop\scriptstyle p=\,{\rm a\,\, critical\,\,point\,\, of}\,\,{\mathcal H} ({\xi})}
\frac{e^{c{\mathcal H}({\xi})(p)}}{\left[{\rm det}\, L_p({\xi})
\right]^{\frac{1}{2}}}
\mbox{.}
\label{DHf}
\end{equation}
In addition the critical points of ${\mathcal H} ({\xi})$ are those where
$d {\mathcal H}({\xi})$ vanishes.
\end{theorem}
Here some $G$-invariant Riemannian metric $\langle\,,\,\rangle$ on $X$ has
been selected, and the square-root in (\ref{DHf}) is as  in (\ref{det3}).

Let us refer to some application of the Duistermaat-Heckman localization formula,
a result when it is specialized, for example, to the symplectic structure on orbit space.
In this connection, the basic example of a Hamiltonian $G$-space is that of an orbit ${{\cO}}$
in the dual space ${\mathfrak g}^*$ of ${\mathfrak g}$ under the co-adjoint action of $G$
on ${\mathfrak g}^*$, which is induced by the adjoint action of $G$ on ${\mathfrak g}$.
In this construction $\sigma$ is chosen as the Kirillov symplectic form on $X={{\cO}}$,
and ${\cH}$ is given by a canonical construction.

Of special interest is the case where ${{\cO}}$ is a quotient $G/T$ of a compact, connected Lie group $G$
modulo a maximal torus $T$. The Duistermaat-Heckman formula leads to the Kirillov integral
formula for irreducible representations of $G$, which has relevance for the geometric quantization
theory. Note that the classical works on hypergeometric series have been initiated from
the physics side, in particular those  works are intimately related to the irreducible representations of the
compact group $G =U(n)$ (see for example \cite{Holman,Milne}).

\subsection{\texorpdfstring{$\mD$}{D}-modules and extension of the Riemann-Roch-Hirzebruch formula}
\label{D}

Let $G_{\mathbb C}$ be a complex algebraic group with Lie algebra $g$. Suppose $G_{\mathbb C}$ acts on a smooth
manifold $X$ and $\cF$ is a sheaf on $X$.
A weak action of $G_{\mathbb C}$ on $\cF$ is an action of $G_{\mathbb C}$ on $\cF$, which extends
the action of $G_{\mathbb C}$ on $X$.

Let $\mD_X$ be the sheaf of differential operators on $X$. In the case when $\cF$ is a ${\mD}_X$-module, more
structure has to be involved. Indeed, if $\varphi$ denote the action of $G_{\mathbb C}$ on $X$, we get a morphism
$ d\varphi: g\rightarrow \Gamma(X, \mD_X)$. $\Gamma(X, \mD_X)$ denotes the global
sections of $\mD_{X}$, which means the ring of differential operators defined on all of $X$.
$\Gamma(X, \mD_X)$ acts on itself by the commutator action $ad$, in such a manner that  we get a map
$ad\cdot d\varphi: g\rightarrow$ End $\Gamma(X, \mD_X)$. It is clear that $\mD_X$ has a weak
action on $G_{\mathbb C}$, say $\gamma$, which yields the map $d\gamma: g \rightarrow$ End $\Gamma(X, \mD_X)$.
Therefore, we can require the map $d\varphi$ to be $G_{\mathbb C}$-equivariant and $ad\cdot d\varphi$ to coincide
with $d\gamma$.

An interesting application of the localization formula is the generalization of the Riemann-Roch-Hirzebruch
integral formula to ${\mD}$-modules.
Such an  statement can be found in \cite{Libine} where it uses the language of
${\mD}$-modules -- sheaves of modules over the sheaf of linear
differential operators. Its flavour can be illustrated by the following example.

Let $G_{\mathbb C}$ be a connected complex algebraic linear
reductive Lie group defined over $\mathbb R$ and
acting algebraically on a smooth complex projective variety $X$.
Suppose that $G \subset G_{\mathbb C}$ is a real Lie subgroup
lying between the group of real points $G_{\mathbb C}(\mathbb R)$
and the identity component. Examine the sheaf of sections
${{\cO}}({E})$ of a $G_{\mathbb C}$-equivariant algebraic line bundle
$({E}, \nabla_{E})$ over a $G_{\mathbb C}$-invariant open algebraic subset
$Y\subset X$ with a $G_{\mathbb C}$-invariant algebraic flat connection
$\nabla_{E}$.
Let $Y_{\mathbb R} \subset X$ be an open $G$-invariant subset
(which may or may not be $G_{\mathbb C}$-invariant).
Consider the cohomology spaces
$  \label{hspaces}
H^\ast(Y_{\mathbb R}, {{\cO}}({E})).
$
The classical Riemann-Roch-Hirzebruch formula computes the index of
${E}$, i.e. the alternating sum
\begin{equation}
\sum_p (-1)^p \dim H^p(Y_{\mathbb R}, {{\cO}}({E}))
\label{Sum}
\end{equation}
with $Y_{\mathbb R} = Y = X$.
For general $Y_{\mathbb R}$ and $Y$ these dimensions can be infinite. However,
further work on this problem allows us to regard the vector spaces
$H^\ast(Y_{\mathbb R}, {{\cO}}({E}))$ as the representations of $G$. Thus, as a substitute for the index,
we can use for the character of the virtual representation the alternating sum (\ref{Sum}).
Note that for finite-dimensional representations the value of the character at the identity
element $e \in G$ equals the dimension of the representation.

\section{Circle actions}
\label{Circle}

First of all, it is convenient to begin with the one-dimensional setting of the
circle actions on a symplectic manifold.

Recall that a symplectic manifold is a smooth manifold $X$ equipped with
a symplectic form ${\sigma}$, that is to say a smooth anti-symmetric two-form
$\sigma \in \Gamma(\bigwedge^2 TX^*)$ on $X$ which is both non-degenerate (thus
$\iota_\xi \sigma(x) \neq 0$, whenever $\xi$ is a vector field that is non-vanishing)
and closed (i.e.  ${d\xi = 0}$). The symplectic manifolds are necessarily even-dimensional
(because odd-dimensional anti-symmetric real matrices automatically have a zero eigenvalue and
are thus degenerate). If $X$ is a ${2n}$-dimensional manifold, the Liouville measure
on that manifold is defined as the volume form ${\sigma^n/n!}$, where, as before,  we use
${\sigma^n}$ to denote the $n$-fold wedge product of $\sigma$ with itself, and we identify
the volume forms with measures.

Given a smooth function $\cH: X \rightarrow {\mathbb R}$ on a symplectic manifold (called the
Hamiltonian), one can associate the Hamiltonian vector field $\xi = \xi_\cH \in \Gamma(TX)$,
defined by requiring that Eq. (\ref{Ham})  holds. From the non-degeneracy of $\sigma$,
we see that $\xi$ vanishes precisely at the critical (or stationary) points of the Hamiltonian $\cH$.
We can exponentiate the Hamiltonian vector field to obtain a one-parameter group
$\varphi(t): X \rightarrow X$ of smooth maps for $t \in {\mathbb R}$:
$
\varphi(t) x := e^{tX} x.
$
This is a smooth action of the additive group ${\mathbb R}$; these maps are
symplectomorphisms and in particular preserve Liouville measure.
One can think of $\varphi$ as a homomorphism from ${\mathbb R}$ to the symplectomorphism
group $Symp(X)$ of $X$.

\subsection{Higher elliptic genera of level \texorpdfstring{$N$}{N} and characteristic \texorpdfstring{$q$}{q}-series}
\label{N}

 Let us consider the generalized elliptic genera for a manifold with an $S^1$-action.
Next, let us introduce  the Ruelle-type spectral function of hyperbolic geometry
${\mathcal R}(s)$ \cite{Bytsenko11}. The function ${\mathcal R}(s)$ is an alternating product
of more complicated factors, each one being  the  so-called Patterson-Selberg zeta-function,
\begin{eqnarray}
\prod_{n=\ell}^{\infty}(1- q^{an+\varepsilon})
& = &\cR(s = (a\ell + \varepsilon)(1-i\varrho(\tau)) + 1-a),
\label{R}
\\
\prod_{n=\ell}^{\infty}(1+ q^{an+\varepsilon})
& = & \cR(s = (a\ell + \varepsilon)(1-i\varrho(\tau)) + 1-a + i\sigma(\tau))\,,
\label{R2}
\end{eqnarray}
where $q\equiv e^{2\pi i\tau}$, $\varrho(\tau) =
{\rm Re}\,\tau/{\rm Im}\,\tau$,
$\sigma(\tau) = (2\,{\rm Im}\,\tau)^{-1}$,
$a$ is a real number, $\varepsilon, b\in {\mathbb C}$, $\ell \in {\mathbb Z}_+$.

Let $N$ be a fixed integer greater than 1, and $x$ be a
variable which is a complex number. Suppose that $\mathfrak h$ is the upper half-plane
of the complex numbers, $\tau \in {\mathfrak h}$.
Let $L=2\pi i({\mathbb Z}\tau+{\mathbb Z})$ be a lattice and $\beta = 2\pi i(\frac{k}{N}\tau +
\frac{l}{N})$, $0\leq k< N, 0\leq l< N$, $\beta \neq 0$. It is usefull to introduce the following
function \cite{Hara}:
\begin{eqnarray}
\Phi(x) & = & (1-e^{-x})\prod_{n=1}^\infty \frac{(1-q^ne^{-x})(1-q^ne^x)}{(1-q^n)^2}
\nonumber \\
& = &
(1-e^{-x})\cdot \frac{\cR(s= (1+ix/2\pi\tau)(1-i\varrho(\tau)))}
{\cR(s=1-i\varrho(\tau))}
\nonumber \\
& \times &
\frac{\cR(s= (1-ix/2\pi\tau)(1-i\varrho(\tau)))}
{\cR(s=1-i\varrho(\tau))}\,.
\label{Phi}
\end{eqnarray}
In fact Eq. (\ref{Phi}) is direct analog of K-theory Euler class \cite{Ando}.
The infinite product in Eq. (\ref{Phi}) is similar to the inverse version of the Witten genus,
which is an even function of the variable $z= i2\pi x$ and defines the  cobordism invariance of
oriented manifolds. Also it admits a {\it Spin} structure and takes values
in ${\mathbb Z}[[q]]$.  Let us take
\begin{equation}
 f(x) := e^{(k/N)x}\Phi(x)\Phi(-\beta)/\Phi(x-\beta),
\end{equation}
where the function $f(x)$ is elliptic with respect to a  sublattice $\widetilde L$ of index $N$
in $L$ \cite{Hirzebruch}.  In the case of a compact  almost complex manifold $X$, the
total Chern class $C(TX)$ of $X$ can be written formally as $C(TX) = \prod_{i=1}^d(1+x_i)$.
In this connection the resulting Chern character takes the form
\begin{eqnarray}
&&
\!\!\!\!\!\!\!\!\!\!
Ch (\bigotimes_{n\geq 1}\cS_{q^n}(TX))^{-1}  = \prod_j\prod_{n=1}^\infty
(1-q^ne^{x_j})(1-q^ne^{-x_j})
\nonumber \\
&&
\!\!\!\!\!\!\!\!\!\!
 =
\prod_j\cR(s=(1-ix_j/2\pi\tau)(1-i\varrho (\tau)))\cdot \cR(s=
(1+ix_j/2\pi\tau)(1-i\varrho (\tau)).
\end{eqnarray}
For a real vector bundle $E$ over $X$ the symmetric powers of $E$ become
$\cS_k(E) = \sum_{k\geq 0}\cS^k(E)t^k$. The elliptic genus of level $N$ can be defined as
$
 \varphi_N(X) = \langle \prod_{j=1}^d\frac{x_j}{f(x_j)}, [X]\rangle.
$
The following assertion is known (see for example \cite{Ando}): if $X$ has complex dimension $m$, then
$\varphi_N(X)$ is a modular form of weight $m$ under the actions of $SL(2, {\mathbb Z})$.

\subsection{Properties of rigidity of elliptic operators}

Recall that the elliptic genus has  originated
in order to find the generating series of rigid elliptic operators.
Let as before $X$ be a smooth compact manifold with an action of a group $G$.
Suppose $\cD$ is an elliptic operator on $X$ (it commutes with  the action). The kernel and the cokernel
of $\cD$ are finite representations of $G$ and the Lefschetz number of $\cD$ at $g\in G$ (a character of $G$)
can be written as follows: $\cL_\mD(g) = tr_gker\,\cD - tr_gcoker\,\cD$. We presume that $\cD$ is
rigid with respect to $G$ if $\cL$ is independent of $g$; in other words we say that a character of $G$
is constant. It is evident that in order to proove the rigidity of an elliptic operator with respect to
a general compact connected Lie group action, we need to analyze its rigidity with respect to $S^1$-action.

\begin{remark}
Let $X$ be a spin manifold, then the rigidity of the Dirac operator with respect to $S^1$-action is the
well known $\widehat{A}$-vanishing theorem \cite{Atiyah2}.
Using the two half-spinor representations of $Spin(2k)$, one can get two associated bundles on $X$
denoted by $\{\triangle^+, \triangle^-\}$. For any real representation $E$ of $G$, the Atiyah-Singer index
theorem gives
\begin{equation}
{\rm Index}\, D\,\otimes E = \int_X {\rm ch} E \,\widehat{A}(X),
\end{equation}
where $D\otimes\,E$ is the twisted Dirac operator;
$D\, \otimes E:\, \Gamma(\triangle^+\, \otimes E)\rightarrow \Gamma(\triangle^-\, \otimes E)$.
In addition, the index is defined to be ${\rm Index}\, D\,\otimes E = {\rm dim\,ker}\, D\,\otimes E
- {\rm dim\,coker}\, D\,\otimes E$.
The conjecture of the rigidity of $D\,\otimes TX$ and  its further proof for the compact homogeneous spin
manifolds have been given in \cite{Witten}.

The rigidity of the signature operator with respect to $S^1$-action on the loop space has been conjectured
by Witten; it is a highly non-trivial result.
Motivated by quantum field theory, Witten has conjectured that the elliptic operators associated with the
appropriate elliptic genera are rigid. Witten also has presented the rigidity conjectures for almost
complex manifolds.

In many physical problems the partition functions (and the Poincar\'{e} polynomials) are linked to
the dimensions of an appropriate homologies for topological spaces and admit  infinite-product representations
for the generating functions and the elliptic genera. In this connection,
the symmetry modular properties of the q-deformed partition functionscan be formulated by the specifics  of the
Ramanujan's summation formula for the bilateral hypergeometric function. We define the general bilateral
hypergeometric function ${}_r\psi_s$ as follows \cite{Gasper}:
\begin{equation}
{}_r\psi_s\left[\begin{matrix}
a_1, a_2, \ldots, a_r \cr
b_1, b_2, \ldots, b_s
\end{matrix}; q,z\right]
= \sum_{n=-\infty}^{\infty}\frac{(a_1, a_2, \ldots, a_r; q)_n}{(b_1, b_2, \ldots; q)_n}
(-1)^{(s-r)n}q^{(s-r)n(n-1)/2}z^n.
\label{psi}
\end{equation}
It is assumed that each term of this series is correctly defined.
\footnote{As an example, consider  ${}_3\psi_2$-summation formula \cite{Gasper}. Let $n$
be a non-negative integer; $q$-analog of the Pfaff-Saalsch\"{u}tz summation formula for ${}_3\psi_2$
is
$$
{}_3\psi_2\left[\begin{matrix}
a, b, q^{-n};q, q \cr
c,abq^{1-n}/c
\end{matrix} \right]
= \frac{(c/a; q)_n(c/b; q)_n}{(c; q)_n(c/ab; q)_n}
$$
}
$(a)_n = a(a+1)\cdots (a+n-1)$. The shifted $q$-factorial has the form:
$(a; q)_n = 1$ for $n=0$ and
\begin{eqnarray}
(a; q)_n &=& \frac{(a;q)_\infty}{(aq^n; q)_\infty} = \prod_{m=0}^\infty [(1-aq^m)/(1-aq^{m+n})]
= (1-a)(1-aq)\cdots (1-aq^{n-1}),
\nonumber \\
&{}& n=1,2,\ldots,
\label{aq}
\end{eqnarray}
where $n$ is a non-negative integer.
\end{remark}
Extensive work in the theory of partition identities shows that the basic hypergeometric series provide
the generating functions for numerous families of identities.
Note that the symmetry properties (or more general the modular properties) for elliptic genera
can be analyzed on the basis of (\ref{psi}) and also functional equations for the spectral
functions of hyperbolic geometry $\cR(s)$ (see for example Eq. (4.22) in \cite{BytsenkoChaichian}).

\section{Twisted equivariant K-theory}
\label{Twist}

{\bf Equivariant K-theory.}
At this point we collect some basic properties of equivariant K-theory.
Recall that (classically) topological K-theory asserts the association of any finite CW complex
$X$ to the category of vector bundles on $X$. Since it is an additive category, where all the exact
sequences split, one can define the Grothendieck group $K(X)$ of this category.  This can be achieved
by taking the free Abelian group on all symbols $[E]$ for $E$ a vector bundle on $X$.
$K(X)$ is a ring because one has a tensor product operation on the category of vector bundles (which
commutes with direct sums). As a result, we can  define a functor $K_G(\cdot)$ from compact $G$-spaces
to Abelian groups, which descends to the homotopy category of equivariant spaces. (It is clear that in
the non-equivariant case, $K_G(X)$ is a commutative ring). With these preliminaries recall the following
definition: The equivariant K-group $K_G(X)$ of $X$ is the Grothendieck group of equivariant
vector bundles on $X$.

Some relations between the equivariant K-theory and the ordinary K-theory have to be mentioned:
({\bf i}) There is a functor from equivariant vector bundles to vector bundles.
({\bf ii}) There is a ring-homomorphism $K_G(X) \rightarrow K(X)$.
Let $H \rightarrow G$ be a morphism of compact Lie groups, then
there is a map $K_G(X) \rightarrow K_H(X)$.
({\bf iii}) Let $X$ be a $G$-space and $X/G$ be an ordinary space, then there is a map
$X \rightarrow X/G$. As a result, there is a map of a vector bundle on $X/G$ to an ordinary
vector bundle on $X$.
({\bf iv}) There is a map $K(X/G) \rightarrow K_G(X)$.
\footnote{Note that if $G$ acts freely on $X$, then
the map $K(X/G) \rightarrow K_G(X)$ is an isomorphism.
}

{\bf Twisted version of equivariant K-theory.}
Let us discuss a twisted version of equivariant K-theory for a global quotient.
Assume that $G$ is a finite group and suppose that for a given class $\alpha\in H^2(G,S^1)$
the compact Lie group extension $1\to S^1\to\widetilde{G}_{\alpha}\to G\to 1$ is fulfilled.
Here the group $\widetilde{G}_{\alpha}$ is associated with the structure of a compact Lie group,
where $S^1\to \widetilde{G}_{\alpha}$ is the inclusion of  a closed subgroup. Let $X$ be a finite
$G$-CW complex.

We say that an $\alpha$-twisted $G$-vector bundle on $X$ is a complex vector bundle $E\to X$, if
$S^1$ acts on the fibers through complex multiplication. In addition, the action extends to an
action of $\widetilde{G}_{\alpha}$ on $E$, which covers a given $G$-action on $X$.
The $\alpha$-twisted $G$-equivariant K-theory of $X$ we denote by $K_G(X, \alpha)$, which is
defined as the Grothendieck group of isomorphism classes of $\alpha$-twisted $G$-bundles over $X$.
\footnote{
See \cite{Wang} for the precise definition.}
In the case $\alpha = 0$, the Lie group extension
corresponds to the split extension $G\times S^1$. In addition, any $G$-vector bundle
can be turned  into a $G\times S^1$-bundle by means of a scalar multiplication on the fibers.
It means that a $G\times S^1$-bundle restricts to a $G$-bundle. As a result, we have
$K^\ast_G(X, \alpha =0) = K^\ast_G(X)$.

For a finite group $G$ and $G$-CW complex $X$, decomposition for K-theory takes the form \cite{Adem}:
\begin{equation}
K^\ast_G(X, \alpha)\otimes {\mathbb C}\cong \bigoplus_{\{(g)\vert g\in G\}}(K^\ast(X^{\langle g\rangle})
\otimes L_g^\alpha)^{Z_G(g)}\,.
\label{dec}
\end{equation}
In Eq. (\ref{dec}), $(g)$ is the conjugacy class of $g\in G$, $Z_G(g)$ and denotes the centralizer of $g$
in $G$, while $L^\alpha_G$ denotes the character for the centralizer $Z_G(g)$.

\subsection{The symmetric product}
\label{Sym}

Let $G$ be a finite group and $X$ be a $G$-space, on which a normal subgroup $A$
acts trivially. The $G$-equivariant K-theory of $X$ can be decomposed as a
direct sum of twisted equivariant K-theories of $X$ parametrized by the orbits of
the conjugation action of $G$ on the irreducible representations of $A$
\cite{Gomez}(theorems 3.2 and 3.4).

Now we  turn to the symmetric group and denote by $G=\mS_n$ the symmetric group on $n$ symbols.
As before, we denote the non-trivial class by $\alpha$.
Using the decomposition formula (theorem 3.4 of \cite{Gomez}), one can calculate
$K^*_{\mS_n}(X^n, \alpha)$,  where the group acts on the
$n$-fold product of a manifold $X$ by permutation of coordinates.
At last, formula for the twisted symmetric products takes the form
\begin{eqnarray}
&&
\sum~q^n\chi (K^*_{\mS_n}(X^n, \alpha)\otimes\mathbb C)
=
\prod_{n>0} (1-q^{2n-1})^{-\chi (X)}
+ \prod_{n>0} (1+q^{2n-1})^{\chi (X)}
\nonumber \\
&&
\times
[1+\frac{1}{2}
\prod_{n>0}(1+q^{2n})^{\chi (X)}
- \frac{1}{2}\prod_{n>0}(1-q^{2n})^{\chi (X)}].
\label{final}
\end{eqnarray}
In Eq. (\ref{final}), $\chi (X)$ is the Euler characteristic of $X$.  In terms of spectral
functions $\cR(s)$, Eq. (\ref{final}) acquires the form
\begin{eqnarray}
&&
\sum~q^n\chi (K^*_{\mS_n}(X^n, \alpha)\otimes\mathbb C)
 =  \cR(s= -i\varrho(\tau))^{-\chi(X)} + \cR(s= -i\varrho(\tau)+i\sigma(\tau))^{\chi(X)}
\nonumber \\
&&
\times
[1+ \frac{1}{2}\cR(s= 1-2i\varrho(\tau)+i\sigma(\tau))^{\chi(X)}
-\frac{1}{2}\cR(s= 1-2i\varrho(\tau))^{\chi(X)}]\,.
\end{eqnarray}

\subsection{The orbifold symmetric product}
\label{Orbifold}

We will analyze the Poincar\'{e} polynomial of {\it orbispace} $[X^n/{\mathfrak S}_n]:=[X\times\cdots
\times X/{\mathfrak S}_n]$ \cite{LUX1}, whose  category  objects are $n$-tuples $(x_1, \cdots, x_n)$ of points in $X$.
In addition, the arrows are elements of the form $(x_1,\cdots, x_n; \sigma)$, where $\sigma \in {\mathfrak S}_n$.
$(x_1,\cdots, x_n; \sigma)$ has  its source as $(x_1,\cdots, x_n)$, while its target as $(x_{\sigma(1)}, \cdots,x_{\sigma(n)})$.
\footnote{
The mentioned category is a groupoid for the inverse of $(x_1,\cdots, x_n; \sigma)$, which is
$(x_1,\cdots, x_n; \sigma^{-1})$. This is the  reason why we can consider  $[X^n/{\mathfrak S}_n]$ as an
orbispace.}

{\bf Poincar\'{e} $q$-polynomials.} Let $X$ be a topological space, denoted by $\nP(X,y)$ its Poincar\'{e} polynomial,
$\nP(X,y) = \sum_jb^j(X)y^j$. Here $b^j(X)$ is the $j$-th Betti number of $X$.
The following formula has been proved by Macdonald \cite{Macdonald}
\begin{equation}
\sum_{n=0}^\infty \nP(X^n/\mS_n, y)q^n = \frac{\prod_j(1+qy^{2j+1})^{b^{2j+1}(X)}}
{\prod_j(1-qy^{2j})^{b^{2j}(X)}}.
\label{M}
\end{equation}

For $y=-1$, Eq.(\ref{M}) reduces to the
formula for the Euler characteristic of the symmetric product
$
\sum_{n=0}^\infty \chi(X^n/\mS_n)q^n = (1-q)^{-\chi(X)}$,
which is valid for the  topological space $X$ whose cohomology $H^j(X; {\rm real})$ is finitely generated
for each $j\geq 0$. Using the $\mS_n$-equivariant K-theory of $X^n$, one can arrive at the following
formula for generation functions \cite{LUX1}:
$
\sum_{n=0}^\infty \chi_{\mS_n}(X^n/\mS_n)q^n = \prod_{j>0}(1-q^j)^{-\chi(X)}.
$
In terms of the spectral functions $\cR(s)$, Eq. (\ref{M}) acquires the form
\begin{equation}
\sum_{n=0}^\infty \nP(X^n/\mS_n, y)q^n = \frac{\cR(s=(1-i({\rm log}(y))(2j+1)/(2\pi\tau))(1-i\varrho(\tau))
+i\sigma(\tau))^{b_{2j+1}(X)}}{\cR(s=(1-i({\rm log}(y))(j/\pi\tau))(1-i\varrho(\tau))^{b_{2j}(X)}}\,.
\end{equation}

The Poincar\'e orbifold polynomial can be defined as follows \cite{LUX1}
\begin{equation}
\nP_{orb}([X/G],y):= \sum y^j
{\rm rank}\,H_{orb}^j([X/G]; {\rm real}) \equiv \sum y^j b^j_{orb}([X/G])\,,
\end{equation}
where $b_{orb}^j$ is the $j$-th orbifold Betti number.

\begin{remark}
For any arbitrary (fixed) positive integers $m_1, \cdots, m_r$ the following equality holds
\begin{equation}
 \sum_{n_1,\cdots, n_r=0}^\infty F(n_1,\ldots,n_r)z_1^{n_1},\ldots, z_r^{n_r}
 =\sum_{N=0}^\infty\,\,\sum_{\stackrel{n_1m_1+\cdots+n_rm_r =N}{n_1,\ldots,n_r\geq 0}}
F(n_1,\ldots,n_r)z_1^{n_1},\ldots, z_r^{n_r}\,.
\label{sum}
\end{equation}
\end{remark}

For the symmetric product, viewed as an orbifold groupoid $[X^n/{\mathfrak S}_n]$, it follows that
\begin{equation}
H_{orb}^\ast([X^n/{\mathfrak S}_n];{\rm real}) \cong \bigoplus_{\sum jn_j= n} \bigotimes_j H^\ast(X^{n_j};
{\rm real})^{\mathfrak{S}_{n_j}}\,.
\label{Horb}
\end{equation}
Using Eq. (\ref{sum}) for calculating the orbifold Poincar\'e polynomial,
we obtain

\begin{eqnarray}
&& \sum_{n=0}^\infty q^n  \nP_{orb}([X^n/{\mathfrak S}_n],y)
 =
\sum_{n=0}^\infty q^n \left(\sum_{\sum jn_j= n} \prod_j \nP(X^{n_j}/\mathfrak{S}_{n_j},y) \right)
\nonumber \\
&& =
\prod_{j>0}\prod_{n \geq 0} \frac{(1 + q^n y^{2j+1})^{b^{2j+1}(X)}}
{(1 - q^n y^{2j})^{b^{2j}(X)}}\,
\nonumber \\
&& =
\prod_{j>0} \frac{[{\cR} (s = -i(2j+1){\rm log}(y)/(2\pi\tau)(1-i\varrho(\tau))+i\sigma(\tau))
]^{b^{2j+1}(X)}}
{[{\cR} (s = -ij{\rm log}(y)/(\pi\tau)(1-i\varrho(\tau))]^{b^{2j}(X)}}\,.
\label{prod}
\end{eqnarray}
For these formulas to be valid, the cohomology of $X$ must be finitely generated at each $n$. The first line
in (\ref{prod}) is similar to the equality (2.32) of \cite{LUX1}. As before, $b^{j}(X)$ is the $j$-th Betti
number of $X$.

\section*{Acknowledgments}

We are grareful to Anca Tureanu for many discussions and several suggestions.
AAB and AEG would like to acknowledge the Conselho Nacional
de Desenvolvimento Cient\'{i}fico e Tecnol\'{o}gico (CNPq, Brazil) and
Coordenac\~{a}o de Aperfei\c{c}amento de Pessoal de N\'{i}vel Superior
(CAPES, Brazil) for financial support.

\end{document}